\documentclass[epsfig,latexsym,amsfonts,twoside]{article}
\usepackage{amssymb,amsmath,amscd}

\def\part#1{\frac{\partial\phantom{#1}}{\partial#1}}
\newtheorem{thm}{Theorem}

\newtheorem{lem}[thm]{Lemma}

\newenvironment{prf}{\begin{trivlist}\item[]{\bf Proof} }%
{\hfill $\Box$ \end{trivlist}}
{\end{trivlist}}
\newenvironment{rmk}{\begin{trivlist}\item[]{\bf Remark} }%
{\end{trivlist}}
\newenvironment{exm}{\begin{trivlist}\item[]{\bf Example} }%
{\end{trivlist}}

\def\Z{\ifmmode{{\mathbb Z}}\else{${\mathbb Z}$}\fi}
\def\Q{\ifmmode{{\mathbb Q}}\else{${\mathbb Q}$}\fi}
\def\C{\ifmmode{{\mathbb C}}\else{${\mathbb C}$}\fi} 
\def\P{\ifmmode{{\mathbb P}}\else{${\mathbb P}$}\fi} 

\def\H{\ifmmode{{\mathrm H}}\else{${\mathrm H}$}\fi} 

\def\B{\ifmmode{{\cal B}}\else{${\cal B}$}\fi} 
\def\E{\ifmmode{{\cal E}}\else{${\cal E}$}\fi} 
\def\F{\ifmmode{{\cal F}}\else{${\cal F}$}\fi} 
\def\K{\ifmmode{{\cal K}}\else{${\cal K}$}\fi} 
\def\L{\ifmmode{{\cal L}}\else{${\cal L}$}\fi} 
\def\M{\ifmmode{{\cal M}}\else{${\cal M}$}\fi} 
\def\N{\ifmmode{{\cal N}}\else{${\cal N}$}\fi} 
\def\O{\ifmmode{{\cal O}}\else{${\cal O}$}\fi} 
\def\U{\ifmmode{{\cal U}}\else{${\cal U}$}\fi}
\def\X{\ifmmode{{\cal X}}\else{${\cal X}$}\fi} 

\def\Br{\ifmmode{{\mathrm{Br}}}\else{${\mathrm{Br}}$}\fi} 
\def\OG{\ifmmode{\widetilde{\cal M}_4}\else{$\widetilde{\cal M}_4$}\fi} 
\def\D{\ifmmode{{\cal D}_{\mathrm{coh}}^b}\else{${{\cal
    D}_{\mathrm{coh}}^b}$}\fi}
\def\Shah{\ifmmode{\amalg\hspace*{-3.5pt}\amalg}\else{$\amalg\hspace*{-3.5pt}\amalg$}\fi}

\begin{document}

\title{Deformations of holomorphic Lagrangian fibrations\footnote{2000 {\em Mathematics Subject Classification.\/} 53C26; 14D06; 14J60.}}
\author{Justin Sawon}
\date{August, 2005}
\maketitle

\begin{abstract}
Let $X\rightarrow\P^n$ be a $2n$-dimensional projective holomorphic 
symplectic manifold admitting a Lagrangian fibration over
$\P^n$. Matsushita proved that the fibration can be deformed in a
codimension one family in the moduli space $\mathrm{Def}(X)$ of
deformations of $X$. We extend his result by proving that if the
Lagrangian fibration admits a section, then there is a codimension two
family of deformations which also preserve the section.
\end{abstract}

\section{Introduction}

Let $X$ be a $2n$-dimensional compact K{\"a}hler manifold. We say $X$ is a 
{\em holomorphic symplectic manifold\/} if it admits a closed two-form 
$\sigma$ of type $(2,0)$ which is non-degenerate in the sense that 
$\sigma^{\wedge n}$ trivializes the canonical bundle $K_X=\Omega^{2n}$. 
Moreover, we call $X$ {\em irreducible\/} if $\H^0(X,\Omega^2)$ is 
one-dimensional and generated by $[\sigma]$. Huybrechts' notes 
in~\cite{ghj02} provide a comprehensive introduction to the standard 
results on irreducible holomorphic symplectic manifolds.

Let us summarize what is known about fibrations on $X$. Suppose we 
have a proper (holomorphic) surjection $f:X\rightarrow B$ onto a complex 
space $B$ such that the general fibre is connected and $0<\mathrm{dim}B<2n$. 
If $X$ is projective and $B$ is a normal variety 
then Matsushita~\cite{matsushita99,matsushita00i} showed that
\begin{itemize}
\item every irreducible component of a fibre of $f$ is a (holomorphic)
Lagrangian submanifold of $X$; in particular, it is $n$-dimensional,
\item the generic fibre is an abelian variety,
\item $B$ is $n$-dimensional and has only $\Q$-factorial log-terminal
singularities,
\item $K_B^*$ is ample and $B$ has Picard number one.
\end{itemize}
Moreover, if $B$ is smooth then it has the same Hodge numbers as 
$\P^n$ (see~\cite{matsushita05}). In particular, for $n=2$ the base 
is isomorphic to $\P^2$ (a similar result was obtain by 
Markushevich~\cite{markushevich96}). In general it is 
expected that $B$ should be isomorphic to $\P^n$, and more evidence 
is provided by a result of Cho, Miyaoka, and 
Shepherd-Barron~\cite{cmsb02} which states that if $f:X\rightarrow B$
admits a section (or even local sections) then $B\cong\P^n$.
Huybrechts~\cite[Proposition 24.8]{ghj02} extended some of these 
results by dropping the projectivity assumption: he showed that if $X$ 
and $B$ are (smooth) K{\"a}hler manifolds then
\begin{itemize}
\item every fibre of $f$ is (holomorphic) Lagrangian,
\item every smooth fibre is a complex torus,
\item $B$ is $n$-dimensional and projective, 
\item $K_B^*$ is ample and its Picard and second Betti numbers both equal one.
\end{itemize}

In this paper we consider deformations of $f:X\rightarrow B$. It is
known that the Kuranishi space $\mathrm{Def}(X)$ of deformations of
$X$ is a smooth complex manifold of dimension $b_2-2$, where $b_2$ is
the second Betti number of $X$. Under the assumptions that $X$ is
projective and $B\cong\P^n$, Matsushita~\cite{matsushita05} proved
that there is a codimension one submanifold of $\mathrm{Def}(X)$
parametrizing deformations of $X$ which are Lagrangian fibrations over
$\P^n$. We include a proof of this result in Section 3 in order to
set-up our notation. Our main result is Theorem~\ref{deformation_s}:
we prove that if $X$ is projective and $f:X\rightarrow B$ admits a
section, then there is a codimension two submanifold of
$\mathrm{Def}(X)$ parametrizing deformations which are Lagrangian
fibrations and admit sections.

Let $X^{\prime}$ be an arbitrary irreducible holomorphic symplectic 
manifold. Can $X^{\prime}$ be deformed to a Lagrangian fibration? The 
answer in general is unknown. However, the question has been answered in
the affirmative for all known examples of irreducible holomorphic 
symplectic manifolds: see Beauville~\cite{beauville99} for the Hilbert 
schemes of points on a K3 surface, Debarre~\cite{debarre99} for the 
generalized Kummer varieties, and Rapagnetta~\cite{rapagnetta04} for
O'Grady's examples.

The author would like to thank Akira Fujiki, Daniel Huybrechts, and
Daisuke Matsushita for useful conversations. Thanks also to Ryushi
Goto for the invitation to Japan where this result was obtained. The
author is support by NSF grant number 0305865.



\section{Some preliminaries}

We collect together some results that we will need later.

\begin{thm}[Matsushita~\cite{matsushita05}]
\label{directimage}
Let $f:X\rightarrow B$ be a Lagrangian fibration between smooth 
projective manifolds. Then
$$R^jf_*\O_X\cong\Omega^j_B.$$
\end{thm}

Let $X$ be a projective irreducible holomorphic symplectic manifold which
admits a Lagrangian fibration $f:X\rightarrow\P^n$ over projective space.
Let $L:=f^*\O(1)$ be the pull-back of the hyperplane line bundle to $X$.
We can use Theorem~\ref{directimage} to calculate the cohomology of 
$L$.

\begin{lem}
\label{sections}
For $(X,L)$ as above
$$h^0(X,L):=\mathrm{dim}\H^0(X,L)=n+1$$
and all higher cohomology of $L$ vanishes.
\end{lem}

\begin{prf}
The Leray spectral sequence gives
$$\H^i(\P^n,R^jf_*L)\Rightarrow\H^{i+j}(X,L).$$
On the left hand side
\begin{eqnarray*}
R^jf_*L & = & R^jf_*(f^*{\cal O}(1)) \\
 & \cong & {\cal O}(1)\otimes R^jf_*\O_X \\
 & \cong & {\cal O}(1)\otimes\Omega^j_{\P^n} \\
\end{eqnarray*}
by the projection formula and Theorem~\ref{directimage}. We can 
therefore use Bott's generalization of the Borel-Weil 
theorem~\cite{es02} to calculate the spectral sequence, finding
$$\H^i(\P^n,\O(1)\otimes\Omega^j_{\P^n})\cong\C^{n+1}$$
when $i=j=0$, and vanishes otherwise. Therefore the sequence 
collapses, and the lemma follows.
\end{prf}

\begin{rmk}
The proof of the lemma relies on two important assumptions: that the 
base of the fibration is $\P^n$ and that $X$ is projective (note that
the projectivity of $X$ is crucial for Theorem~\ref{directimage}).
\end{rmk}

We need to know how the cohomology of a line bundle behaves in a
family.

\begin{thm}[Grauert and Remmert~\cite{gr84}, page 210]
\label{upper}
Let ${\cal X}\rightarrow\Delta$ be a deformation of a complex manifold
$X={\cal X}_0$, where $\Delta$ is a polydisc of some dimension around 
$0$. Let $V$ be a vector bundle on $\cal X$, with restriction $V_t$ to 
each ${\cal X}_t$. Then for $i\geq 0$,
$$h^i({\cal X}_t,V_t):=\mathrm{dim}\H^i({\cal X}_t,V_t)$$
is an upper semi-continuous function of $t\in\Delta$. In other words
$$A_{i,d}:=\{t\in\Delta|h^i({\cal X}_t,V_t)\geq d\}$$
is a closed analytic subset of $\Delta$ for each integer $d$.
\end{thm}



\section{Deforming fibrations}

Let $X$ be an irreducible holomorphic symplectic manifold of dimension 
$2n$. For the following statements, see Section 22 of Huybrechts' notes
in~\cite{ghj02}. Denote by
$${\cal X}\rightarrow(\mathrm{Def}(X),0)$$
the Kuranishi family parametrizing local deformations of $X={\cal X}_0$.
We think of $(\mathrm{Def}(X),0)$ as the germ of a complex space. It is 
smooth (deformations are unobstructed) of dimension $b_2-2$, where $b_2$ 
is the second Betti number of $X$. Note that when we deform $X$ as a 
complex manifold it remains holomorphic symplectic and irreducible; for 
small deformations it also stays K{\"a}hler.

The germ $(\mathrm{Def}(X),0)$ can be represented by a contractible open
set, and therefore for each $t\in\mathrm{Def}(X)$ we can choose an 
isomorphism
$$\psi_t:\H^2({\cal X}_t,\Z)\rightarrow\H^2(X,\Z)$$
known as a {\em marking}. Let
$$Q_X:=\{[\alpha]\in\P(\H^2(X,\C))|q_X(\alpha)=0\mbox{ and }q_X(\alpha+\bar{\alpha})>0\}$$
be the {\em period domain\/} of $X$, where $q_X$ is the 
Beauville-Bogomolov quadratic form on $\H^2(X,\Z)\otimes\C$. Then by the 
Local Torelli Theorem (see Beauville~\cite{beauville83}), the {\em period 
map\/}
$${\cal P}_X:(\mathrm{Def}(X),0)\rightarrow (Q_X,[\sigma])$$
which takes $t$ to $[(\psi_t\otimes\C)(\sigma_t)]$ is a local isomorphism.

In this section we prove the following result.

\begin{thm}[Matsushita~\cite{matsushita05}, Corollary 1.7]
\label{deformation}
Let $X$ be an irreducible projective holomorphic symplectic manifold 
which admits a Lagrangian fibration $f:X\rightarrow\P^n$ over 
projective space. There is a codimension one submanifold 
$\Delta^f\subset\mathrm{Def}(X)$ (containing zero) which parametrizes 
holomorphic symplectic manifolds $X^{\prime}$ which admit Lagrangian 
fibrations over $\P^n$.
\end{thm}

\begin{prf}
Let $L:=f^*\O(1)$ be the pull-back of the hyperplane line bundle to 
$X$. Then $f$ is clearly given by the morphism
$$\phi_L:X\rightarrow\P(\H^0(X,L)^*)$$
induced by the linear system of $L$. Our aim is to extend $L$ and this 
map to certain deformations of $X$.

Let $c_1:=c_1(L)\in\H^2(X,\Z)\cap\H^{1,1}(X)$ be the first Chern class 
of $L$, and let
$$\Delta^f:=\{t\in\mathrm{Def}(X)|\psi_t^{-1}(c_1)\in\H^2({\cal X}_t,\Z)\mbox{ is of type }(1,1)\}.$$
If $t\in\Delta^f$, then $\psi_t^{-1}(c_1)$ is orthogonal to $\sigma_t$
with respect to $q_{{\cal X}_t}$, or equivalently, $c_1$ is orthogonal
to $(\psi_t\otimes\C)(\sigma_t)$ with respect to $q_X$. Therefore 
$\Delta^f$ maps isomorphically to a neighbourhood of $[\sigma]$ in
$$Q_X^f:=\{[\alpha]\in Q_X|q_X(\alpha,c_1)=0\}$$
which is of codimension one in $Q_X$. Hence $\Delta^f$ must be 
codimension one in $\mathrm{Def}(X)$.

The exponential exact sequence on ${\cal X}_t$ gives
$$\ldots\rightarrow 0\rightarrow\H^1({\cal X}_t,\O^*)\rightarrow\H^2({\cal X}_t,\Z)\stackrel{\beta}{\rightarrow}\H^2({\cal X}_t,\O)\rightarrow\ldots$$
since $\H^1({\cal X}_t,\O)$ vanishes. For $t\in\Delta^f$, 
$\psi_t^{-1}(c_1)$ is of type $(1,1)$ and so in the kernel of $\beta$. 
It therefore comes from a unique holomorphic line bundle 
$\L_t\in\H^1({\cal X}_t,\O^*)$. Moreover, we can always choose a 
representative $(U,0)$ of the germ $(\mathrm{Def}(X),0)$ such that
$U\cap\Delta^f$ is connected and simply-connected, so there exists a 
line bundle $\L$ over ${\cal X}|_{\Delta^f}$, restricting to $\L_t$ on 
each fibre ${\cal X}_t$.

By Theorem~\ref{upper}, $h^i({\cal X}_t,\L_t)$ is an upper 
semi-continuous function of $t\in\Delta^f$. So for $t$ in a 
neighbourhood of zero
$$h^i({\cal X}_t,\L_t)\leq h^i({\cal X}_0,\L_0).$$
But the right hand side is $h^i(X,L)$, which by Lemma~\ref{sections}
is $n+1$ for $i=0$ and zero otherwise. Therefore
$$h^i({\cal X}_t,\L_t)=0$$
for $i>0$. On the other hand, the Euler characteristic
$$\chi({\cal X}_t,\L_t):=\sum_{i=0}^{2n}(-1)^ih^i({\cal X}_t,\L_t)$$
is constant (it is given by the Hirzebruch-Riemann-Roch formula) 
and equal to $\chi(X,L)=n+1$. Therefore
$$h^0({\cal X}_t,\L_t)=n+1$$
for all $t$ in a neighbourhood of zero. The linear system of $\L_t$
therefore gives a map
$$\phi_t:{\cal X}_t\rightarrow\P(\H^0({\cal X}_t,\L_t)^*)\cong\P^n.$$

A priori $\phi_t$ is only a rational map, so we must show that the 
linear system of $\L_t$ has no base-points. This follows
since the specialization
$$\phi_L=f:X\rightarrow\P(\H^0(X,L)^*)\cong\P^n$$
is a morphism and thus has no base-points. More specifically, we can 
extend a basis $\{s_0,\ldots,s_n\}$ of sections of $L$ to a set of 
$n+1$ sections $\{S_0,\ldots,S_n\}$ of $\L$ over 
${\cal X}|_{\Delta^f}\rightarrow\Delta^f$. The locus
$$C:=\{t\in\Delta^f|S_0|_{{\cal X}_t},\ldots,S_n|_{{\cal X}_t}\mbox{ simultaneously vanish at some point of }{\cal X}_t\}$$
is a closed subset of $\Delta^f$ which does not contain $0$. So for 
$t\in\Delta^f$ in a neighbourhood of zero, the sections 
$S_0|_{{\cal X}_t},\ldots,S_n|_{{\cal X}_t}$ generate the fibres of $\L_t$.

We have proved that 
$$\phi_t:{\cal X}_t\rightarrow\P^n$$
is a morphism, and it is therefore a Lagrangian fibration by the 
results of Matsushita (and Huybrechts) cited in the introduction.
\end{prf}

\section{Fibrations with sections}

We will use the following result (note that in Voisin's paper $X$ is
the Lagrangian submanifold of $Y$).

\begin{thm}[Voisin~\cite{voisin92}]
\label{voisin}
Let $X$ be a holomorphic symplectic manifold and $g:Y\hookrightarrow X$ 
a Lagrangian submanifold. Let $g_t^*$ be the composition
$$g^*\circ(\psi_t\otimes\C):\H^2({\cal X}_t,\C)\rightarrow\H^2(Y,\C)$$
where $\psi_t$ is a marking of ${\cal X}_t$. Then the inclusion
$Y\hookrightarrow X$ deforms to a Lagrangian submanifold 
${\cal Y}_t\rightarrow{\cal X}_t$ if and only if $g_t^*(\sigma_t)=0$.

We can rephrase this as follows. Let $L_Y\subset\H^2(X,\C)$ be the 
orthogonal complement of $\mathrm{ker}g^*$ with respect to the 
Beauville-Bogomolov quadratic form $q_X$. Since $\H^{2,0}(X)\subset\mathrm{ker}g^*$ 
and $L_Y$ can be defined over $\Q$, $L_Y$ must be of type $(1,1)$. Then 
$Y\hookrightarrow X$ deforms to ${\cal Y}_t\rightarrow{\cal X}_t$ if and 
only if $L_Y$ is preserved under the deformation, i.e.\ if and only if
$$t\in\Delta:=\{t\in\mathrm{Def}(X)|\psi^{-1}_t(L_Y)\in\H^2({\cal X}_t,\C)\mbox{ is of type }(1,1)\}.$$
\end{thm}

We wish to look at deformations of a Lagrangian fibration 
$f:X\rightarrow B$ which admits a section $s$. Note that the existence
of the section implies that $B\cong\P^n$~\cite{cmsb02}. However, we 
will keep the projectivity of $X$ as a hypothesis. Let $Y$ denote the 
image of $s$, and $g:Y\hookrightarrow X$ the inclusion. Since $Y$ is 
isomorphic to $B\cong\P^n$, it must be holomorphic Lagrangian, as 
$$\sigma|_Y\in\H^0(Y,\Omega^2)\cong\H^{2,0}(\P^n)=0.$$

\begin{thm}
\label{deformation_s}
Let $X$ be an irreducible projective holomorphic symplectic manifold 
which admits a Lagrangian fibration $f:X\rightarrow B$ and a global 
section $s$. There is a codimension two submanifold 
$\Delta^{fs}\subset\mathrm{Def}(X)$ (containing zero) which 
parametrizes holomorphic symplectic manifolds $X^{\prime}$ which admit 
Lagrangian fibrations and global sections.
\end{thm}

\begin{prf}
Firstly, the existence of the section implies $B\cong\P^n$, and 
so by Theorem~\ref{deformation} there exists a codimension one family 
$\Delta^f\subset\mathrm{Def}(X)$ over which the fibration deforms.

By Theorem~\ref{voisin} we also have
$$\Delta^s:=\{t\in\mathrm{Def}(X)|\psi^{-1}_t(L_Y)\in\H^2({\cal X}_t,\C)\mbox{ is of type }(1,1)\}$$
over which the Lagrangian submanifold $Y\hookrightarrow X$ deforms.
Now $\Delta^s$ maps isomorphically to a neighbourhood of $[\sigma]$ in
$$Q_X^s:=\{[\alpha]\in Q_X|q_X(\alpha,L_Y)=0\}=\{[\alpha]\in Q_X|\alpha\in\mathrm{ker}g^*\}.$$
Since $Y\cong\P^n$, $\H^2(Y,\C)$ and $L_Y$ are one-dimensional, so 
$Q_X^s$ and $\Delta^s$ are codimension one in $Q_X$ and 
$\mathrm{Def}(X)$ respectively.

Let $\Delta^{fs}:=\Delta^f\cap\Delta^s$. It maps isomorphically to a 
neighbourhood of $[\sigma]$ in
$$Q_X^{fs}:=Q_X^f\cap Q_X^s=\{[\alpha]\in Q_X|q_X(\alpha,c_1)=0\mbox{ and }\alpha\in\mathrm{ker}g^*\}.$$
Observe that $L:=f^*\O(1)$ must satisfy
$$c_1(L)^{n+1}=0\in\H^{2n+2}(X,\Z),$$
which implies $q_X(c_1)=0$
(see~\cite{ghj02}) and thus $[c_1]\in Q_X^f$. On the other hand
$$g^*(c_1)=c_1(\O(1))\neq 0,$$
so $c_1\not\in\mathrm{ker}g^*$ and
$[c_1]\not\in Q_X^s$. This suffices to show that $Q_X^f$ and $Q_X^s$
intersect transversely, and thus $Q_X^{fs}$ and $\Delta^{fs}$ are
codimension two in $Q_X$ and $\mathrm{Def}(X)$ respectively.

Finally observe that if $t\in\Delta^{fs}$ then ${\cal X}_t$ is a 
Lagrangian fibration over $\P^n$ and it contains a deformation 
${\cal Y}_t$ of $Y$. By specialization, ${\cal Y}_t$ is a section of 
the fibration. More precisely, the set of values $t$ for which 
${\cal Y}_t$ maps isomorphically to $\P^n$ is open in $\Delta^{fs}$ 
and contains zero. This completes the proof.
\end{prf}

\begin{exm}
Let $S$ be a K3 surface which contains a smooth genus $g$ curve $C$, 
$g\geq 2$, such that $\mathrm{Pic}S$ is generated by $[C]$. Let ${\cal
  C}\rightarrow |C|\cong\P^g$ be the family of curves linearly
equivalent to $C$. Then the compactified Jacobian
$$J^0:=\overline{\mathrm{Jac}}^0({\cal C}/\P^g)$$
is a deformation of the Hilbert scheme of $g$ points on 
$S$~\cite{beauville99}. Moreover, $J^0$ is a Lagrangian fibration 
over $\P^g$ which admits a section. Observe that we can deform $S$ in
a 19-dimensional family, while keeping the genus $g$ curve: these produce
deformations of $J^0$ keeping the Lagrangian fibration and section. On the
other hand, the space of deformations of the Hilbert scheme of points on 
$S$ is 21-dimensional. This is in agreement with 
Theorem~\ref{deformation_s}.

We find similar agreement for Lagrangian fibrations which are deformations
of the generalized Kummer varieties (see Debarre~\cite{debarre99}, or the 
author's paper~\cite{sawon03}).
\end{exm}

\begin{rmk}
Let $f:X\rightarrow B$ be a Lagrangian fibration. If the fibres of $X$ 
are all reduced and irreducible then there exists the (Altman-Kleiman)
compactified relative Picard scheme
$$P:=\overline{\mathrm{Pic}}^0(X/B).$$
If the fibres of $P$ are reduced and irreducible we also have 
$$X^0:=\overline{\mathrm{Pic}}^0(P/B).$$
Both $P$ and $X^0$ admit global sections. Let $U$ and $U^0$ denote the
open subsets of $X$ and $X^0$, respectively, given by removing all
singular fibres; then $U$ is a torsor over $U^0$. In fact in cases
where the singular fibres of $X$ are not too complicated, $X$ itself
is a torsor over $X^0$; moreover $P$ and $X^0$ are holomorphic
symplectic manifolds (see~\cite{sawon05}). In~\cite{sawon04} it was
shown that torsors $X$ over $X^0$ are parametrized by the
one-dimensional space $\H^2(P,\O^*)$ of gerbes on $P$, with $X^0$
being the unique fibration which admits a global section. This agrees
with Theorem~\ref{deformation} and Theorem~\ref{deformation_s}, which
together imply that fibrations which admit sections must be
codimension one inside the family of all fibrations.
\end{rmk}

\begin{flushleft}
Department of Mathematics\hfill sawon@math.sunysb.edu\\
SUNY at Stony Brook\hfill www.math.sunysb.edu/$\sim$sawon\\
Stony Brook NY 11794-3651\\
USA\\
\end{flushleft}


\begin{thebibliography}{XXX}

\bibitem{beauville83} A. Beauville,
{\em Vari{\'e}t{\'e}s K{\"a}hl{\'e}riennes dont le 1{\'e}re classe de
Chern est nulle\/},
Jour. Diff. Geom. {\bf 18} (1983), 755--782.

\bibitem{beauville99} A. Beauville,
{\em Counting rational curves on K3 surfaces\/},
Duke Math. J. {\bf 97} (1999), no. 1, 99--108.


\bibitem{cmsb02} K. Cho, Y. Miyaoka, and N. Shepherd-Barron,
{\em Characterizations of projective space and applications to complex
  symplectic manifolds\/},
in Higher Dimensional Birational Geometry, Advanced Studies in Pure
Mathematics {\bf 35}, 2002, 1--88.

\bibitem{debarre99} O. Debarre,
{\em On the Euler characteristic of generalized Kummer varieties\/},
Amer. J. Math. {\bf 121} (1999), no. 3, 577--586.

\bibitem{es02} M. Eastwood and J. Sawon,
{\em The Borel-Weil theorem for complex projective space},
in Invitations to Geometry and Topology, Oxford University Press
(2002), 126--145.

\bibitem{gr84} H. Grauert and R. Remmert,
{\em Coherent analytic sheaves}, 
Grundlehren der Mathematischen Wissenschaften {\bf 265}, 
Springer-Verlag, Berlin, 1984.


\bibitem{ghj02} M. Gross, D. Huybrechts, and D. Joyce,
{\em Calabi-Yau manifolds and related geometries\/},
Springer Universitext, 2002.





\bibitem{markushevich96} D. Markushevich,
{\em Lagrangian families of Jacobians of genus 2 curves\/},
J. Math. Sci. {\bf 82} (1996), no. 1, 3268--3284.

\bibitem{matsushita99} D. Matsushita,
{\em On fibre space structures of a projective irreducible symplectic
manifold},
Topology {\bf 38} (1999), No. 1, 79--83. 
Addendum, Topology {\bf 40} (2001), no. 2, 431--432.

\bibitem{matsushita00i} D. Matsushita,
{\em Equidimensionality of Lagrangian fibrations on holomorphic 
symplectic manifolds},
Math. Res. Lett. {\bf 7} (2000), no. 4, 389--391.


\bibitem{matsushita05} D. Matsushita,
{\em Higher direct images of dualizing sheaves of Lagrangian 
fibrations},
Amer. J. Math. {\bf 127} (2005), no. 2, 243--259.




\bibitem{rapagnetta04} A. Rapagnetta,
{\em Topological invariants of O'Grady's six dimensional irreducible 
symplectic variety\/},
preprint {\bf math.AG/0406026}.

\bibitem{sawon03} J. Sawon,
{\em Abelian fibred holomorphic symplectic manifolds\/},
Turkish Jour. Math. {\bf 27} (2003), no. 1, 197--230. 

\bibitem{sawon04} J. Sawon,
{\em Derived equivalence of holomorphic symplectic manifolds\/},
in Algebraic Structures and Moduli Spaces, CRM Proc. Lecture Notes
{\bf 38}, AMS, 2004, 193--211.

\bibitem{sawon05} J. Sawon,
{\em Twisted Fourier-Mukai transforms for holomorphic symplectic
  four-folds\/}, 
preprint.

\bibitem{voisin92} C. Voisin,
{\em Sur la stabilit{\'e} des sous-vari{\'e}t{\'e}s lagrangiennes 
des vari{\'e}t{\'e}s symplectiques holomorphes\/},
Complex projective geometry, 294--303, LMS Lecture Note Ser. {\bf 179},
Cambridge Univ. Press, 1992.





\end{thebibliography}
\end{document}